\documentclass[11pt,twoside]{article}
\usepackage[mathscr]{euscript}
\usepackage{secdot}
\usepackage[utf8]{inputenc}
\usepackage{amsmath, amsthm, amscd, amsfonts, amssymb, color}
\usepackage[bookmarksnumbered,plainpages]{hyperref}

\date{}

\begin{document}

\centerline{}

\centerline {\Large{\bf A STUDY ON DELTA COMPACT SPACES}}
\centerline{}
\centerline{\textbf{Sanjay Roy${^*}$}}
\centerline{Department of Mathematics, Uluberia College}
\centerline{Uluberia, Howrah- 711315,  West Bengal, India}
\centerline{e-mail: sanjaypuremath@gmail.com}
\centerline{}
\centerline{\textbf{Srabani Mondal}}
\centerline{Department of Mathematics, Uluberia College}
\centerline{Uluberia, Howrah- 711315,  West Bengal, India}
\centerline{e-mail: srabanimondal531@gmail.com}
\centerline{}
\centerline{\textbf{Shrobana Sinha Roy}}
\centerline{Department of Mathematics, Uluberia College}
\centerline{Uluberia, Howrah- 711315,  West Bengal, India}
\centerline{e-mail: sinharoy.shrobana@gmail.com}
\centerline{}

\centerline{\textbf{Bobi Mandal}}
\centerline{Department of Mathematics, Uluberia College}
\centerline{Uluberia, Howrah- 711315,  West Bengal, India}
\centerline{e-mail: bobimandal767@gmail.com}
\newcommand{\mvec}[1]{\mbox{\bfseries\itshape #1}}
\centerline{}
\newtheorem{Theorem}{\quad Theorem}[section]

\newtheorem{definition}[Theorem]{\quad Definition}

\newtheorem{theorem}[Theorem]{\quad Theorem}

\newtheorem{remark}[Theorem]{\quad Remark}

\newtheorem{corollary}[Theorem]{\quad Corollary}

\newtheorem{note}[Theorem]{\quad Note}

\newtheorem{lemma}[Theorem]{\quad Lemma}

\newtheorem{example}[Theorem]{\quad Example}
\newtheorem{notation}[Theorem]{\quad Notation}

\newtheorem{result}[Theorem]{\quad Result}
\newtheorem{conclusion}[Theorem]{\quad Conclusion}

\newtheorem{proposition}[Theorem]{\quad Proposition}
\newtheorem{prop}[Theorem]{\quad Property}

\begin{abstract}
\textbf{\emph{The aim of this paper is to introduce the concept of Delta-Compact spaces along with some basic properties of it. Here, we try to\\
establish the behavior of Delta-Compact spaces under the continuous mapping. Finally, we define another concept namely, Locally Delta-Compactness. }}
\end{abstract}
{\bf Keywords:}  \emph{Delta Compact Set, Delta Closed mapping, Locally Delta compact Set.}\\
\textbf{2010 Mathematics Subject Classification:} 54D30, 54D45, 54F99\\ 
\\${^*}$Corresponding Author

\section{\textbf{Introduction}}
If we take the interior of closure of an open set, formerly known as regular open set, how it can be associated to a given topological space. To find out this association, in 1980, T. Noiri \cite{Noiri} introduced the concept of $\delta$-closed set with the help of $\delta$-closure and defined a $\delta$-open set as a complement of $\delta$-closed set. There he defined the concept of $\delta$-continuous function and established it's few properties. In 1993, S. Raychaudhuri and M. N. Mukherjee \cite{Raychaudhuri} generalized the concepts of $\delta$-continuity and $\delta$-open set and called these general definitions as $\delta$-almost continuity and $\delta$-preopen set respectively. In 1993, N. Palaniappan and K. C. Rao \cite{Palaniappan} introduced a concepts of regular generalized closed set and regular generalized open set. Thereafter many works have been done on another concepts of $\delta$-semiopen set \cite{Park, Caldas 1}, $\delta$-set \cite{Saleh}, Semi$^{*}\delta$-open set \cite{Missier 2}, regular$^{*}$ open set \cite{Missier 3}.

In this paper, we define the concepts of $\delta$-open set and $\delta$-closed set in another way and later we have shown that these two concepts are equivalent with the definitions which were defined by T. Noiri \cite{Noiri}. Compactness is an important property of a topological space. Here we have established the concept of Delta compact space. As every compact space is a Delta compact space but not vice versa, each properties which holds in Delta compact space, holds in compact space but not vice versa. So, it is a more advanced topic than compact Space. \\
In this paper, we also introduce the notion of locally Delta compactness along with some properties and examples with the help of regular open sets.

\section{Preliminaries}
Here we have given a few basic definitions and theorems which will be needed in the sequel.

\begin{definition}\cite{Noiri}
An open set $U$ of a topological space $(X,\tau)$ is said to be a \textbf{Regular Open Set} if $U=$ $int[cl(U)]$.
\end{definition}

\begin{definition}\cite{Noiri}
A subset $V$ of a topological space $(X,\tau)$ is said to be \textbf{regular closed} if $V=$ $cl[int(V)]$.
\end{definition}

\begin{theorem}\label{th2}
$V$ is a regular closed set iff $X\setminus V$ is a regular open set.
\end{theorem}

\begin{proof}
Let us suppose that $V$ is a regular closed set i.e, $V=$ $cl[int(V)]$.\\
We have to show that, $X\setminus V$ is a regular open set i.e, $X\setminus V=$ $int[cl(X\setminus V)]$.\\
Let $x\in X\setminus V$. Then $x\notin V= cl[int V]$. So, there exists a nbd $U_x$ of $x$ such that $U_x\cap int V=\phi\\
\Rightarrow U_x\cap V=\phi\\
\Rightarrow U_x\subseteq X\setminus V \\
\Rightarrow U_x\subseteq cl( X\setminus V )\\
\Rightarrow x\in int[cl( X\setminus V )]$. Thus $X\setminus V \subseteq int[cl( X\setminus V )]$. Again let $x\in int[cl(X\setminus V)]$. Then there exists a nbd $U_x$ of $x$ such that $U_x\subseteq cl(X\setminus V)=X\setminus int V\Rightarrow U_x\cap int V=\phi\Rightarrow x\notin cl[int V]=V\Rightarrow x\in X\setminus V.$ Thus $ int[cl( X\setminus V )]\subseteq X\setminus V.$

Thus, we get $X\setminus V=$ $int[cl(X\setminus V)]$.\\
Hence, $X\setminus V$ is a regular open set.
\end{proof}

\begin{result}\label{R1}
Let $(X, \tau)$ be a topological space and $U\in \tau$. Then $int[cl(U)]$ is regular open.
\end{result}

\begin{proof}
Let $int[cl(U)] = V$. We have to show that, $V =int[cl(V)]$.\\
Since $int[cl(U)] = V$, $V\subseteq cl(U)$ i.e., $cl(V) \subseteq cl(U)$ i.e., $int[cl(V)] \subseteq int[cl(U)]=V$.\\
Again, $V\subseteq$ $cl(V)\Rightarrow$ $int(V)\subseteq$ $int[cl(V)]$.
As, $V$ is open, $int(V)= V$. So, $V\subseteq$ $int[cl(V)]$.
Thus we get $V=int[cl(V)]$.
So, $V= int[cl(U)]$ is a regular open set.
\end{proof}

\begin{definition}\cite{Noiri}
A point $x$ of a topological space $X$ is called the
$\delta$-cluster  point of $A\subseteq X$  if $A \cap int [Cl(U)]\neq\phi$ for every open set $U$ of $X$ containing
$x$.\\
The set of all $\delta$-cluster points of $A$  is called the $\delta$-closure of $A$ and is denoted by
$Cl_{\delta} A$.\\
A subset $A$ of $X$ is called $\delta$-closed if $A=Cl_{\delta} A$. 
The complement of a $\delta$-closed set is called  $\delta$-open set.
\end{definition}

\section{Delta Open ($\delta$-Open) Set and Delta Closed ($\delta$-Closed) Set}

\begin{definition}
A subset $U$ of $(X,\tau)$ is called \textbf{$\delta$-open set} if for each $x\in U$, $\exists$ a regular open set $P$ such that $x\in P$ $\subseteq U$.
\end{definition}

\begin{theorem}
Let $A$ be a $\delta$-open set in a topological space $(X, \tau)$. Then $A$ is a $\delta$-open set according to the definition of T. Noiri \cite{Noiri}.
\end{theorem}

\begin{proof}
Let $A$ be a $\delta$-open set in a topological space $(X, \tau)$ according to our definition. We now show that $X\setminus A= cl_{\delta}(X\setminus A)$. Let $x\in cl_{\delta}(X\setminus A)$. Then $(X\setminus A)\cap int[cl(U)]\neq\phi$ for every open set $U$ containing $x$. If $x\in A$ then there exists a regular open set $V$ such that $x\in V\subseteq A$. So, $X\setminus A\subseteq X\setminus V$. Now $(X\setminus A)\cap int[cl(V)]\subseteq (X\setminus V)\cap int[cl(V)]= (X\setminus V)\cap V=\phi$ which contradicts our assumetion that $x\in cl_{\delta}(X\setminus A)$. So $x\in X\setminus A$. Thus $cl_{\delta}(X\setminus A)\subseteq (X\setminus A)$. Obviously, $(X\setminus A)\subseteq cl_{\delta}(X\setminus A)$. Hence $A$ is a $\delta$-open set according to the definition of T. Noiri \cite{Noiri}.\\
Conversely, let $A$ be a $\delta$-open set according to the definition of T. Noiri \cite{Noiri}. Then $X\setminus A= cl_{\delta}(X\setminus A)$. Let $x\in A$. Then $x\notin X\setminus A$i.e., $x\notin cl_{\delta}(X\setminus A)$. So, there exists an open set $U$ containing $x$ such that $(X\setminus A)\cap int[cl(U)]=\phi$ i.e., $int[cl(U)]\subseteq A$. Now  $int[cl(U)]$ is regular open set containing $x$ by the Result \ref{R1}. Thus for each $x\in A$, there exists a regular open set $V$ such that $x\in V\subseteq A$. Hence $A$ is $\delta$-open according to our definition.
\end{proof}

\begin{example}
If $V$ be a regular open set in $(X,\tau)$, then $V\cap Y$ may or may not be a regular open set in $(Y, \tau_{Y})$, where $Y$ is a subset of $X$.

Let us take the usual topology on real numbers $\mathbb{R}$.\\
 Let $V = (1,2)$ be open in $\mathbb{R}$ and $Y$ = $[1,1.5]$ a subspace of $\mathbb{R}$. Then $V\cap Y$ = $(1, 1.5]$. So,  $cl_{Y}(V\cap Y)$ = $[1, 1.5]$ and $int_{Y}cl_{Y}(V\cap Y)$ = $[1, 1.5] \neq$ $V\cap Y$.
\end{example}

\begin{theorem}\label{th1}
If $V$ be a regular open set in $(X,\tau)$ and $Y$ be an open set in $(X,\tau)$ then $V\cap Y$ is also a regular open set in $(Y,\tau_{Y})$.
\end{theorem}

\begin{proof}
Since $V$ is a regular open set in $(X,\tau)$, $V$ = $int(cl V)$. Then $V\in \tau$. So,$V\cap Y \in$ $\tau_{Y}$. Thus, $V\cap Y \subseteq$ $int_{Y}cl_{Y}(V\cap Y)$. Again let $x\in$ $int_{Y}cl_{Y}(V\cap Y)$. Then there exists $W\in \tau$ such that $x\in$ $W\cap Y$ $\subseteq cl_{Y}(V\cap Y)$ $\subseteq cl(V\cap Y)$ $\subseteq cl(V)$. Since $W\in \tau$ and $Y\in \tau$, $W\cap Y\in \tau$. So, $x\in$ $int(cl V)$ = $V$. Again $x\in Y$, thus, $x\in$ $V\cap Y$. Hence $V\cap Y$ = $int_{Y}cl_{Y}(V\cap Y)$.
\end{proof}

\begin{theorem}
If $U$ be a $\delta$- open set in $(X,\tau)$ and $Y\in$ $\tau$, then $U\cap Y$ is a $\delta$- open set in $(Y,\tau_{Y})$.
\end{theorem}

\begin{proof}
Let $U$ be a $\delta$- open set in $(X,\tau)$. We now show that $U\cap Y$ is also a $\delta$- open set in $(Y,\tau_{Y})$. Let $x\in$ $U\cap Y$. Since $x\in$ $U\cap Y$, then $x\in U$ and $x\in Y$. Again since $x\in U$ and $U$ is a $\delta$- open set, then by definition, $\exists$ a regular open set $P$ in $(X,\tau)$ such that $x\in P$ $\subseteq U$. Then $x\in$ $P\cap Y \subseteq$ $U\cap Y$.Since $P$ is regular open set in $(X,\tau)$ and $Y\in \tau$, then by the previous theorem we an say that $P\cap Y$ is a regular open set in $(Y,\tau_{Y})$. Hence $U\cap Y$ is a $\delta$- open set in $(Y,\tau_{Y})$.
\end{proof}

\begin{definition}
A subset $V$ of $(X,\tau)$ is called $\delta$-closed set if $V$ is the intersection of all regular closed set containing $V$.
\end{definition}

\begin{theorem}
$V$ is $\delta$-open set iff $X\setminus V$ is $\delta$-closed set.
\end{theorem}

\begin{proof}
Suppose $V$ is a $\delta$-open set. Since $V$ is $\delta$-open, then for each $x\in V$, $\exists$ a regular open set $U$ such that $x\in U$ $\subseteq V$ and so $X\setminus U$ $\supseteq X\setminus V$.\\
Again since $U$ is regular open set, so by the Theorem \ref{th2}, $X\setminus U$ is regular closed set. Therefore there exists a regular closed set $X\setminus U$ containing $X\setminus V$ and also if $x\notin X\setminus V$ i.e., if $x\in V$, there exists a regular closed set $X\setminus U$ containing $X\setminus V$ such that $x\notin X\setminus U$. So, $X\setminus V$ is the intersection of all regular closed set containing $X\setminus V$ Hence, $X\setminus V$ is $\delta$-closed.\\
Conversely, let us suppose that $X\setminus V$ is a $\delta$-closed set.\\
Let $x\in V$. Then $x\notin X\setminus V$. Since $X\setminus V$ is $\delta$-closed, there exists a regular closed set $F$ containing $X\setminus V$ such that $x\notin F$. So, there exists a regular open set $X\setminus F$ such that $x\in X\setminus F\subseteq V$. Hence $V$ is $\delta$-open.
\end{proof}

\begin{note}
Let $A$ be a $\delta$-closed set in a topological space $(X, \tau)$. Then $A$ is a $\delta$-closed set according to the definition of T. Noiri \cite{Noiri}.
\end{note}

\section{Delta Compact($\delta$-Compact) Spaces}
\begin{definition}
Let $(X,\tau)$ be a topological space. $X$ is called a \textbf{$\delta$-Compact Space} if every open cover of $X$ by $\delta$-open subset of $X$ has a finite subcover.
\end{definition}

\begin{definition}
Let $(X,\tau)$ be a topological space and $Y\subseteq X$. Then $Y$ is called a $\delta$-compact if every $\delta$-open cover of $(Y,\tau_{Y})$ has a finite subcover.
\end{definition}

\begin{theorem}
An open subset $Y$ of a topological space $(X,\tau)$ is $\delta$-compact if and only if every $\delta$-open cover of $Y$ by the $\delta$-open subset of $X$ has a finite subcover.
\end{theorem}

\begin{proof}
Let $Y$ be a $\delta$- compact subset of $(X,\tau)$. Let $\lbrace G_{\alpha}: \alpha\in \Lambda\rbrace$ be a $\delta$- open cover of $Y$, where each $G_{\alpha}$ is $\delta$- open set in $(X,\tau)$ $\forall \alpha\in \Lambda$.\\
Then, $Y\subseteq \cup_{\alpha\in\Lambda}G_{\alpha}$ $\Rightarrow Y\subseteq \cup_{\alpha\in\Lambda}G_{\alpha}$ $\cap Y$, where each $G_{\alpha}\cap$ $Y$ is a $\delta$- open set in $(Y,\tau_{Y})$ by the Theorem   \ref{th1}.\\
Therefore, by $\delta$- compactness of $Y$, $\exists$ a finite subcollection $\Lambda_{0}$ of $\Lambda$ such that $Y\subseteq \cup_{\alpha\in\Lambda_{0}}G_{\alpha}$ $\cap Y$ $\Rightarrow Y\subseteq$ $\cup_{\alpha\in\Lambda_{0}}G_{\alpha}$.\\
Thus, if $Y$ is $\delta$- compact then every $\delta$- open cover of $Y$ by the $\delta$- open set of $(X,\tau)$ has a finite subcover.\\
Conversely, let $\lbrace Y_{\alpha}: \alpha\in \Lambda\rbrace$ be a $\delta$- open cover of $Y$ by the $\delta$- open sets of $\tau_{Y}$.Therefore, $Y\subseteq$ $\cup_{\alpha\in \Lambda}$ $Y_{\alpha}$.\\
Since $Y_{\alpha}$ is a $\delta$- open set in $\tau_{Y}$ and $Y$ is open, $Y_{\alpha}$ is a $\delta$- open set in $\tau$ $\forall \alpha\in \Lambda$.\\
Therefore, $\lbrace Y_{\alpha}: \alpha\in \Lambda\rbrace$ is a $\delta$- open cover of $Y$ by the $\delta$- open sets of $\tau$.\\
Then by the given condition, $\exists$ a finite subcover $\Lambda_{0}$ of $\Lambda$ such that $Y\subseteq$ $\cup_{\alpha\in \Lambda}$ $Y_{\alpha}$.\\
So by the definition of $\delta$- compact space, $Y$ is $\delta$- compact.\\
Hence, the theorem is done.
\end{proof}

\begin{theorem}
A topological space $(X,\tau)$ is $\delta$- compact iff for every collection $\lbrace F_{\alpha}: \alpha\in \Lambda\rbrace$ of $\delta$-closed sets in $X$ possessing finite intersection property $($ f.i.p $)$, $\cap_{\alpha\in \Lambda}$ $F_{\alpha}$ of the entire collection is non-empty.
\end{theorem}

\begin{proof}
Let $(X,\tau)$ be a $\delta$- compact space and $\lbrace F_{\alpha}: \alpha\in \Lambda\rbrace$ be a collection of $\delta$- closed sets possessing f.i.p. We now show that $\cap_{\alpha\in \Lambda}$ $F_{\alpha}$ $\neq \phi$.\\
 If possible let $\cap_{\alpha\in \Lambda}$ $F_{\alpha}$ = $\phi$\\
 Then $X$ = $X\setminus  \cap_{\alpha\in \Lambda}$ $F_{\alpha}$ = $\cup_{\alpha\in\Lambda}(X\setminus F_{\alpha})$, a $\delta$- open cover of $X$. Since $X$ is $\delta$- compact, $\exists$ a finite subcollection $\Lambda_{0}$ of $\Lambda$ such that $X$ = $\cup_{\alpha\in\Lambda_{0}}(X\setminus F_{\alpha})$ = $X\setminus \cap_{\alpha\in\Lambda_{0}}F_{\alpha}$.\\
 Therefore, $\cap_{\alpha\in\Lambda_{0}}F_{\alpha}$ $\neq\phi$, which is a contradiction as $\lbrace F_{\alpha}: \alpha\in \Lambda\rbrace$ satisfies f.i.p.\\
 Hence, $\cap_{\alpha\in\Lambda}F_{\alpha}$ $\neq\phi$.\\
 Conversely, let us assume that the given condition holds.
We now show that $X$ is $\delta$- compact.\\
For this, let $\lbrace G_{\alpha}: \alpha\in \Lambda\rbrace$ be a $\delta$- open cover of $X$.\\
Then $X$ = $\cup_{\alpha\in\Lambda}G_{\alpha}$ $\Rightarrow$ $X \setminus \cup_{\alpha\in\Lambda}G_{\alpha}$ = $\phi$ $\Rightarrow$ $\cap_{\alpha\in \Lambda}(X \setminus G_{\alpha})$ $= \phi$.\\
Then $\exists$ a finite subcollection $\Lambda_{0}$ of $\Lambda$ such that $\cap_{\alpha\in \Lambda_{0}}(X \setminus G_{\alpha})$ $= \phi$ that is, $X \setminus \cup_{\alpha\in\Lambda_{0}}G_{\alpha}$ = $\phi$ that is, $X= \cup_{\alpha\in\Lambda_{0}}G_{\alpha}$.\\
Thus the space $(X,\tau)$ is $\delta$- compact.
\end{proof}

\begin{theorem}\label{th7}
A $\delta$-closed subset of a $\delta$-compact space is $\delta$-compact.
\end{theorem}

\begin{proof}
Let $(X,\tau)$ be a $\delta$-compact topological space and let $Y$ be a $\delta$-closed subset of $X$.\\
We now show that, $Y$ is $\delta$-compact.\\
Let, $\lbrace G_{\alpha}: \alpha\in \Lambda\rbrace$ be a $\delta$-open cover of $Y$, where each $G_{\alpha}$ is $\delta$-open set in $(X,\tau)$ $\forall \alpha\in \Lambda$.
Then $Y\subseteq$ $\cup_{\alpha\in \Lambda}G_{\alpha}$.
So, $X\subseteq$ $(X\setminus Y)\cup$ $(\cup_{\alpha\in \Lambda}G_{\alpha})$.
Since $X$ is $\delta$-compact, $\exists$ a finite collection $\Lambda_{0}$ of $\Lambda$ such that $X\subseteq$ $(X\setminus Y)\cup$ $(\cup_{\alpha\in \Lambda_{0}}G_{\alpha})$ and so  $Y\subseteq$ $\cup_{\alpha\in \Lambda_{0}}G_{\alpha}$.
Hence every $\delta$-open cover $\lbrace G_{\alpha}: \alpha\in \Lambda\rbrace$ of $Y$ has a finite subcover.
So, $Y$ is $\delta$-compact.
This completes the proof.
\end{proof}

\begin{theorem}\label{th8}
In a regular space every closed set is a $\delta$-closed set.
\end{theorem}

\begin{proof}
Let $(X,\tau)$ be a regular space and $A\subseteq$ $X$ be closed.\\
We have to show that $A$ is $\delta$-closed i.e, $X\setminus A$ is $\delta$-open.\\
Let $x\notin A$. Then $x\in X\setminus A$.\\
As $A$ is closed, $X\setminus A$ is open.\\
Since the space is regular, $\exists$ an open set $U$ such that $x\in U \subseteq$ $cl(U) \subseteq$ $X\setminus A$.\\
Now, $U\subseteq$ $int[cl(U)]\subseteq$ $cl(U)$. So $x\in$ $int[cl(U)]\subseteq$ $X\setminus A$.\\ Therefore $X\setminus A$ is $\delta$-open set as $int[cl(U)]$ is a regular open set by the Result \ref{R1}.\\
Hence $A$ is $\delta$-closed. 
This completes the proof.
\end{proof}

\begin{corollary} \label{cor1}
In a regular space every open set is also a $\delta$-open set.
\end{corollary}

\begin{theorem}\label{th4}
In a regular $T_{2}$ space, for any two points $x$ and $y$, $\exists$ two disjoint $\delta$-open sets $U$ and $ V$ such that $x\in U$ and $y\in V$.
\end{theorem}

\begin{proof}
Let $(X,\tau)$ be a regular $T_{2}$ space. Let $x$ and $y$ be any two points of $X$.\\
Since the space is $T_{2}$, $\exists$ two disjoint open sets $U$ and  $V$ such that $x\in U$ and $y\in V$.\\
As $x\in U$, $x\notin X\setminus U$ and as $y\in V$,  $y\notin X\setminus V$.\\
Since  $X\setminus U$ and $X\setminus V$ are closed sets
and the space is regular,  $X\setminus U$ and  $X\setminus V$ are $\delta$-closed sets.
Then $U$ and $V$ are $\delta$-open.
So for any two points $x$ and $y$, $\exists$ two disjoint $\delta$-open sets $U$ and $V$ such that $x\in U$ and $y\in V$.\\
\end{proof}

\begin{theorem}\label{th3}
The intersection of finitely many regular open sets in any topological spaces is regular open.
\end{theorem}

\begin{proof}
Let $ (X,\tau) $ be any topological spaces. $ P_1, P_2, P_3 ,\ldots,P_{m}$ be $ m $ regular open sets in $ X $.\\
We take $P$ = $ P_{1} $ $ \cap $  $ P_{2} $ $ \cap $  $ P_{3} $ $ \cap $ $\ldots$ $ \cap $ $ P_{m} $.\\
We now show that $ P $ is regularly open in $ X $, i.e, $ \cap_{i=1}^{m}P_{i} $= int[cl($ \cap_{i=1}^{m}P_{i} $)].\\
Let $ x $ $ \in $ int[cl( $ \cap_{i=1}^{m}P_{i} $)].
Then $ \exists $ $ U\in \tau $ such that, $x\in U \subseteq $ cl($ \cap_{i=1}^{m}P_{i} )$. Then $x\in U \subseteq$ cl($P_{_i}$) for all  $i= 1,2,\ldots,m.$ So $x\in int[cl(P_{i})]= P_i$ for all  $i= 1,2,\ldots, m$ as each $P_{i}$ is regularly open in $X$. So $x\in\cap_{i=1}^{m}P_{i}$.\\
Thus, int[cl( $ \cap_{i=1}^{m}P_{i} $)] $\subseteq$ $ \cap_{i=1}^{m}P_{i} $.\\
Again since $ \cap_{i=1}^{m}P_{i} \subseteq $ cl($\cap_{i=1}^{m}P_{i}$), 
$ \cap_{i=1}^{m}P_{i}\subseteq int[cl(  \cap_{i=1}^{m}P_{i} )]$ ($\because$ each $P_{i}$ is open).\\
Hence, int[cl( $ \cap_{i=1}^{m}P_{i} $)] = $ \cap_{i=1}^{m}P_{i} $.\\
Consequently, $P$ is regular open in $(X,\tau)$.
\end{proof}

\begin{theorem}\label{th5}
The arbitrary union of $ \delta$-open sets in any topological spaces is $ \delta$-open.
\end{theorem}

\begin{proof}
Let $ (X,\tau) $ be any topological spaces and $ \{G_\alpha: \alpha\in \Lambda\}$, a family of $\delta$-open sets in $ X $.
Let $G= \cup_{\alpha\in \Lambda} G_{\alpha}$ and $x\in G$.
Then $x$ belongs to atleast one of the sets, say $G_{\beta}$ where $\beta\in \Lambda$.
Since $G_{\beta}$ is $\delta$-open and $x\in G_{\beta}$, $\exists$ a regular open set $P$ such that $x\in P \subseteq G_{\beta}$.
Therefore $P \subseteq G_{\beta}\subseteq G$.
Since $x$ is arbitrary and $x\in P\subseteq G$, $G$ is $\delta$-open in $(X,\tau)$.
\end{proof}

\begin{theorem}\label{th6}
The intersection of finitely many $\delta$-open sets in any topological spaces is $\delta$-open.
\end{theorem}

\begin{proof}
Let $ (X,\tau) $ be any topological spaces. $ G_1, G_2, G_3 ,\ldots,G_{m}$ be $m$ $\delta$-open sets in $ X $.\\
Let $G$= $\cap_{i=1}^{m} G_{i}$ and $x\in G$.\\
Then, $x\in G_{i}$, for each $i= 1,2,\ldots,m$.\\
Since $G_{1}$ is $\delta$-open, $\exists$ a regular open set $P_{1}$ such that $x\in P_{1}\subseteq G_{1}$.\\
Again, since $G_{2}$ is $\delta$-open, $\exists$ a regular open set $P_{2}$ such that $x\in P_{2}\subseteq G_{2}$.\\
$\ldots\ldots$\\
Continuing in a similar way, since $G_{m}$ is $\delta$-open, $\exists$ a regular open set $P_{m}$ such that $x\in P_{m}\subseteq G_{m}$.\\
Let $P$= $ P_{1} $ $ \cap $  $ P_{2} $ $ \cap $  $ P_{3} $ $ \cap $ $\ldots$ $ \cap $ $ P_{m} $.
 Then by the Theorem \ref{th3}, $P$ is reguarly open in $X$.\\
Now $x\in P\subseteq P_{i}\subseteq G_{i}$ for all $i=1, 2, \cdots m$.
Then $x\in P\subseteq G_{1} \cap G_{2} \cap G_{3} \cap \ldots \cap  G_{m}$.\\
That is,  $x\in P\subseteq G$.
Hence $G$ is $\delta$-open in $(X,\tau)$.
\end{proof}

\begin{theorem}\label{th9}
Let $(X,\tau)$ be any $T_{3}$ topological spaces.Then if $B$ is a $\delta$-compact subset of $X$ and $x$ be any point in $X$ such that $ x $ can be strongly separeted from every point $ y $ of $B$, then $ x $ and $ B $ can also be strongly separated in $ (X,\tau) $ by the $\delta$-open sets.
\end{theorem}

\begin{proof}
If $ B = \phi$, the proof is trivial.\\
If $ B\neq \phi $, let $ y\in B $.\\
Then $ x$ and $y $ are strongly seperated in $ (X,\tau) $.\\
So by the Theorem \ref{th4}, $\exists$ two $\delta$-open sets $ U_{y}$and $V_{y} $ in $ (X,\tau) $ such that $y\in U_{y}$ and $x\in V_{y}$ with  $ U_{y} \cap V_{y} = \phi$.\\
Now $ \lbrace U_{y}: y\in B\rbrace $ is an open cover of $B$. Then $ \lbrace U_{y}: y\in B\rbrace $ is an $\delta$-open cover of $B$ by the Corollary \ref{cor1}\\
Since $ B $ is $ \delta$-compact, $\exists$ a finite subcollection $ \lbrace U_{y_{1}}, U_{y_{2}}, \ldots, U_{y_{n}}\rbrace $ of $ \lbrace U_{y}: y\in B\rbrace $, that also covers $B$.
Now, corresponding to this subcollection, the $ \delta$-open  sets containing $x$ are 
$ V_{{y_{1}}}, V_{{y_{2}}},\ldots,V_{{y_{n}}}$. \\
Let, $ U$=$ \cup_{i=1}^{n}U_{y_{i}}$ and $ V $= $\cap_{i=1}^{n}V_{y_{i}}$.\\
Then by the Theorems \ref{th5} and \ref{th6}, both $ U$ and $ V $ are $\delta$-open sets in $(X,\tau)$ such that $x\in V$ and $B\subseteq U$ and $U\cap V = \phi$.\\
Thus $x$ and $B$ are strongly seperated in $(X,\tau)$ by the $\delta$-open sets $U$ and $V$.
\end{proof}

\begin{theorem}
In every $\delta$- compact $T_{3}$ space, for any two disjoint closed sets $A$ and $B$ $\exists$ two disjoint $\delta$- open sets $U$ and $V$ such that $A\subseteq$ $U$ and $B\subseteq$ $V$.
\end{theorem}

\begin{proof}
Let $(X,\tau)$ be a $\delta$- compact  $T_{3}$ space. Let $A$ and $B$ be two disjoint closed sets. Since $A$ and $B$ are closed subsets of regular space $X$,by the Theorem \ref{th8} $A$ and $B$ are $\delta$-closed sets. Therefore by the Theorem \ref{th7} $A$ and $B$ are $\delta$- compact.  \\
Let $x\in A$. Since $A$ and $B$ are disjoint, then $x\notin B$.\\
 Since the space $(X,\tau)$ is  $T_{3}$ with $x\notin B$, then by the Theorem \ref{th9} we can say that $\exists$ two disjoint $\delta$- open sets $U_{x}$ and $V_{x}$ such that $x\in$ $U_{x}$ and $B\subseteq$ $V_{x}$.\\
Now $A\subseteq\cup_{x\in A} U_{x}$. Since $A$ is $\delta$- compact, there exist finite sub collection $\{U_{x_i}: i=1, 2, \cdots,n\}$  of $\delta$- open sets such that  $A\subseteq\cup_{i=1}^n U_{x_i}$. Let $U=\cup_{i=1}^n U_{x_i}$ and $V=\cap_{i=1}^n V_{x_i}$
Then $U\cap$ $V=$ $\emptyset$ and by the Theorems \ref{th5} and \ref{th6}  $U$ and $V$  are $\delta$- open sets.\\
Hence in a $\delta$- compact, $T_{3}$ space for any two disjoint closed sets $A$ and $B$ $\exists$ two disjoint $\delta$- open sets $U$ and $V$ such that $A\subseteq$ $U$ and $B\subseteq$ $V$.
\end{proof}

\begin{corollary}
Every $\delta$- compact $T_{3}$ space is $T_{4}$.
\end{corollary}

\begin{theorem}\label{th15}
Every $\delta$-compact subset of a $T_{3}$ space is $\delta$-closed.
\end{theorem}

\begin{proof}
Let $(X,\tau)$ be a $T_{3}$ space and $Y$ be a $\delta$-compact subset of $X$.\\
If $Y= X$, there is nothing to prove.
If $Y\neq X$, then let $x\in X\setminus Y$.\\
Since $X$ is $T_{3}$, i.e, regular and $T_{2}$, $x$ can be strongly separated from each point of $Y$. Thus, $x$ can be strongly separated from $Y$ itself in $(X,\tau)$.\\
So, $\exists$ two $\delta$-open sets $U$ and $V$ in $(X,\tau)$ such that $x\in U$, $Y\subseteq V$ and $U\cap V = \phi$. So $x$ is not a limit point in $Y$.
Thus $Y$ is closed in $(X,\tau)$.\\
Hence $Y$ is $\delta$-closed in $(X,\tau)$ as $X$ is regular.
\end{proof}

\begin{corollary}
In a $\delta$-compact $T_{3}$ space, $X$, a set $Y(\subseteq X)$ is $\delta$-compact iff $Y$ is $\delta$-closed, i.e, in such spaces, the concept of $\delta$-closedness and  $\delta$-compactness for subsets coincide.
\end{corollary}

\begin{theorem}
The union of finite collection of $\delta$-compact subsets of a topological space is $\delta$-compact.
\end{theorem}

\begin{proof}
Let $ (X,\tau) $ be any topological spaces and $ Y_1, Y_2, Y_3 ,\ldots,Y_{m}$, finitely many $\delta$-compact subsets in $ X $.\\
Let $Y$= $\cup_{i=1}^{m} Y_{i}$.\\
We have to show that $Y$ is $ \delta$-compact.\\
Let $G$ be a cover of $Y$ by sets which are $\delta$-open in $(X,\tau)$.\\
Then each $Y_{i}$ is covered by $G$ for $i=1,2,\ldots,m$.\\
$\because$ Each $Y_{i}$ is $\delta$-compact, $\exists$ a finitely many subcover $G_{i}$ of $G$ for each $Y_{i}$ $i=1,2,\ldots,m$.\\
Then $G_{1}\cup G_{2}\cup G_{3}\cup\ldots\cup G_{m}$ is a finite subcollection of $G$ which also covers $Y$.\\
Hence $Y$ is $\delta$-compact.
\end{proof}

\begin{theorem}\label{th10}
The intersection of finite collection of $\delta$-closed sets is $\delta$-closed.
\end{theorem}

\begin{proof}
Let $ (X,\tau) $ be any topological spaces and $ Y_1, Y_2, Y_3 ,\ldots,Y_{m}$, finitely many $\delta$-closed sets in $ (X,\tau) $.\\
Let $Y= \cap_{i=1}^{m} Y_{i}$.\\
We have to show that $Y= \cap_{i=1}^{m} Y_{i}$ is $ \delta$-closed.
i.e., $X\setminus Y= X\setminus \cap_{i=1}^{m} Y_{i}$ is $\delta$-open.\\
Since each $Y_{i}$, for $i=1,2,\ldots,m$ is $\delta$-closed in $(X,\tau)$, each
$X\setminus Y_{i}$ for $i=1,2,\ldots,m$ is $\delta$-open in $(X,\tau)$.\\
Then $\cup_{i=1}^{m}(X\setminus Y_{i})$ is also $\delta$-open in $(X,\tau)$\;$( \because$ arbitrary union of $\delta$-open sets is $\delta$-open $)$.\\
Thus, $X\setminus \cap_{i=1}^{m}Y_{i}$ is $\delta$-open in $(X,\tau)$.\\
This shows that, $\cap_{i=1}^{m}Y_{i}$ is $\delta$-closed in $(X,\tau)$.
\end{proof}

\begin{theorem} \label{th11}
The intersection of finite collection of $\delta$-compact subsets of a $T_{3}$ topological space is $\delta$-compact.
\end{theorem}

\begin{proof}
Let $ (X,\tau) $ be any $T_{3}$ topological spaces and $ Y_1, Y_2, Y_3 ,\ldots,Y_{m}$, finitely many $\delta$-compact subsets in $ (X,\tau) $.\\
Let $Y= \cap_{i=1}^{m} Y_{i}$.\\
We now show that $Y$ is $\delta$-compact in $(X,\tau)$.\\
Since each $Y_{i}$ is $\delta$-compact and $X$ is $T_{3}$ and every $\delta$-compact subsets of a $T_{3}$ space is $\delta$-closed,\\
We can say, each $Y_{i}$ is $\delta$-closed.\\
Hence, by the Theorem \ref{th10}, we can say that $Y$ is $\delta$-closed in $(X,\tau)$.
\end{proof}

\begin{theorem}
If $A$ is closed and $B$ is $\delta$-compact in a $\delta$-compact $T_{3}$ space $X$, then $A\cap B$ is $\delta$-compact.
\end{theorem}

\begin{proof}
Since the space $X$ is $T_{3}$, i.e, regular and $T_{2}$, the set $A$ is $\delta$-closed.\\
Also, $X$ is $\delta$-compact, implies $A$ is $\delta$-compact.\\
Thus, by the Theorem \ref{th11}, $A\cap B$ is $\delta$-compact in $X$.
\end{proof}

\section{{$\delta$}-Compactness and Continuity }
In this section, we will try to examine the behaviour of $\delta$-compact spaces under continuous mapping.
\begin{example}
Let, $f:\mathbb{R}^{2}\rightarrow \mathbb{R}^{2}$ be defined by $f(x)= x^{2}$ , $\forall x\in \mathbb{R}$.\\
Now, $f^{-1}(0,1)$= $(-1,0)\cup (0,1)$. \\
Clearly, $(0,1)$ is regularly open.\\
But, cl$[(-1,0)\cup (0,1)]$= $[-1,1]$ and\\
int$[cl((-1,0)\cup (0,1))]$= $(-1,1)\neq (-1,0)\cup (0,1) $\\
This shows that $(-1,0)\cup (0,1)$ is not regularly open.
\end{example}

\begin{note}
If $f:X\rightarrow Y$ is a continuous mapping and $U$ is a regular open set in $Y$, then $f^{-1}$($U$) may not be regularly open set in $X$.
\end{note}

\begin{theorem}\label{th12}
Let $f$ be an open, continuous mapping from $(X,\tau)$ to $(Y,\sigma)$ and $U$ be a regular open set in $Y$. Then $f^{-1}(U)$ is regularly open in $X$.
\end{theorem}

\begin{proof}
Let $f:(X,\tau)\rightarrow (Y,\sigma)$ be an open, continuous mapping.\\
We now show that, $f^{-1}(U)$ is regularly open in $X$.
i.e., int[cl($f^{-1}(U)$)]= $f^{-1}(U)$.\\ 
Let $x\in$ int[cl($f^{-1}(U)$)].
Then $\exists$ $W\in \tau$ such that\\
 $x\in W\subseteq$ cl($f^{-1}(U)$)$\Rightarrow$ $ f(x)\in f(W)\subseteq$ $f(cl(f^{-1}(U))$ $\subseteq cl(f(f^{-1}(U)))$ $\Rightarrow$ $f(x)\in f(W)\subseteq cl(U)$.\\
Since $f$ is an open mapping, $f(W)$ is open in $(Y,\sigma)$.\\
So, $f(x)\in int[cl(U)] \Rightarrow$ $ f(x)\in U $($\because U$ is regularly open in Y)$\Rightarrow x\in f^{-1}(U).$\\
Thus $int[cl(f^{-1}(U))]\subseteq$ $f^{-1}(U)$.\\
Again we know $f^{-1}(U)\subseteq cl(f^{-1}(U))\Rightarrow$ $int[f^{-1}(U)]\subseteq int[cl(f^{-1}(U))]$\\
$\Rightarrow$ $f^{-1}(U)\subseteq int[cl(f^{-1}(U))]$.\\
Hence we have $f^{-1}(U)= int[cl(f^{-1}(U))]$.\\
Consequently, $f^{-1}(U)$ is regularly open in $(X,\tau)$.
\end{proof}

\begin{theorem}\label{th13}
Let $f$ be an open, continuous mapping from $(X,\tau)$ to $(Y,\sigma)$ and $U$ be a $\delta$-open set in $Y$. Then $f^{-1}(U)$ is $\delta$-open in $X$.
\end{theorem}

\begin{proof}
Let $f:(X,\tau)\rightarrow (Y,\sigma)$ be an open, continuous mapping.\\
We now show that, $f^{-1}(U)$ is $\delta$-open in $X$.
i.e., for every $x\in f^{-1}(U)$, there exists a regular open set $V$ in $X$ such that $x\in V\subseteq f^{-1}(U)$.\\ 
Let $x\in f^{-1}(U)$. Then $f(x)\in U$
Since $U$ is $\delta$- open, there exists a regular open set $W$ in $Y$ such that $f(x)\in W\subseteq U$. So, $x\in f^{-1}(W)\subseteq f^{-1}(U).$ Now since $W$ is regular open in $Y$ and $f$ is open continuous mapping, by the Theorem \ref{th12} $f^{-1}(W)$ is regular open in $X$. Hence $f^{-1}(U)$ is $\delta$-open in $X$.
\end{proof}

\begin{theorem}\label{th14}
Let $f$ be an open, continuous mapping from a $\delta$-compact space $(X,\tau)$ to a topological space $(Y,\sigma)$. Then $f(X)$ is $\delta$-compact in $(Y,\sigma)$.
\end{theorem}

\begin{proof}
Let $\lbrace U_{\alpha}:\alpha\in \Lambda\rbrace$ be a $\delta$-open cover of $f(X)$ in $(Y,\sigma)$. Then $f(X)\subseteq \cup_{\alpha\in \Lambda}U_\alpha$, i.e., $X\subseteq f^{-1}(\cup_{\alpha\in \Lambda}U_\alpha)=\cup_{\alpha\in \Lambda}f^{-1}(U_\alpha)$,\\
Since $f$ is open and continuous, by the previous Theorem \ref{th13}, each $f^{-1}(U_{\alpha})$ is $\delta$-open in $X$.\\
Thus $\lbrace f^{-1}(U_{\alpha}):\alpha\in \Lambda\rbrace$ is a $\delta$-open cover of $X$.\\
Since $X$ is $\delta$-compact, then $\exists$ a finite subcollection $\lbrace f^{-1}(U_{\alpha_{1}}), f^{-1}(U_{\alpha_{2}}),\ldots,f^{-1}(U_{\alpha_{m}})\rbrace$ of $\lbrace f^{-1}(U_{\alpha}):\alpha\in \Lambda\rbrace$ which also covers $X$.\\
Thus $X\subseteq$ $\cup_{i=1}^{m}f^{-1}(U_{\alpha_{i}})$.\\
This implies, $f(X)\subseteq$ $\cup_{i=1}^{m}U_{\alpha_{i}}$.\\
Therefore, $\lbrace U_{\alpha_{i}}:i=1,2,\ldots,m\rbrace$ is a finite subcollection of $\lbrace U_{\alpha}:\alpha\in \Lambda\rbrace$ which covers $f(X)$.\\
Hence, $f(X)$ is $\delta$-compact in $(Y,\sigma)$.
\end{proof}

\begin{definition}
A mapping $f:(X,\tau)\rightarrow (Y,\sigma)$ is said to be a \textbf{$\delta$-closed mapping} if it maps from a $\delta$-closed subset of the domain set to a $\delta$-closed subset of the co-domain set.
\end{definition}

\begin{theorem}\label{th16}
Any open continuous mapping from a $\delta$-compact space to any $T_{3}$ topological space is a $\delta$-closed mapping.
\end{theorem}

\begin{proof}
Let $(X,\tau)$ be a $\delta$-compact space and $(Y,\sigma)$ be any $T_{3}$ topological space.\\
Also,let $f:(X,\tau)\rightarrow (Y,\sigma)$ be a open,continuous mapping.\\
We now show that, $f$ is a $\delta$-closed mapping.\\
Let $F$ be a $\delta$-closed subset of $X$. By the Theorem \ref{th7}, it is clear that $F$ is $\delta$-compact in $X$.\\
So, by the Theorem \ref{th14}, $f(F)$ is $\delta$-compact in $(Y,\sigma)$.\\
Also, the space $Y$ is $T_{3}$ and every $\delta$-compact subset in a $T_{3}$ space is $\delta$-closed by the Theorem \ref{th15}. So $f(F)$ is $\delta$-closed in $(Y,\sigma)$.\\
This concludes that $f$ is a $\delta$-closed mapping. 
\end{proof}

\begin{theorem}
Any open continuous mapping from a $\delta$-compact $T_{3}$ space to any $T_{3}$ topological space is a closed mapping.
\end{theorem}

\begin{proof}
Let $(X,\tau)$ be a $\delta$-compact $T_{3}$ space and $(Y,\sigma)$ be any $T_{3}$ topological space.\\
Also,let $f:(X,\tau)\rightarrow (Y,\sigma)$ be a open,continuous mapping.\\
We now show that, $f$ is a closed mapping.
Let $F$ be a closed subset of $X$. Since $X$ is $T_{3}$, i.e, regular and $T_{2}$, $F$ is $\delta$-closed in $X$.
So, by the procedure of the Theorem \ref{th16} it clearly shows that $f(F)$ is $\delta$-closed in $(Y,\sigma)$, that is $f(F)$ is closed in $(Y,\sigma)$.\\
Hence $f$ is a closed mapping. 
\end{proof}

\section{Locally {$\delta$}-Compactness }

\begin{definition}
A topological space $(X,\tau)$ is said to be \textbf{locally $\delta$-compact} at $x\in X$ if there exists a neighbourhood of $x$ which is $\delta$-compact in $X$. If $(X,\tau)$ is locally $\delta$-compact at every point of $ X$ then X is called locally $\delta$-compact
\end{definition}

\begin{example}
Every $\delta$-compact space is a locally $\delta$-compact space. In fact,
let $(X,\tau)$ be a $\delta$-compact space and $x\in X$.
Then $X$ is itself a $\delta$-compact neighbourhood of $x$ in $(X,\tau)$.
Therefore $X$ is locally $\delta$-compact at $x$.
Consequently, $(X,\tau)$ is locally $\delta$-compact.
\end{example}

\begin{example}
Any discrete space $(X,\tau)$ is locally $\delta$-compact. Let $(X,\tau)$ be a discrete space and $x\in X$.
Then $\lbrace x\rbrace$ is a $\delta$-compact neighbourhood of $x$ in $(X,\tau)$.\\
Therefore $X$ is locally $\delta$-compact at $x$.
Consequently, $(X,\tau)$ is locally $\delta$-compact.
\end{example}

\begin{theorem}
Every closed subspace of a locally $\delta$-compact $T_{3}$ space is locally $\delta$-compact.
\end{theorem}

\begin{proof}
Let $Y$ be a closed subspace of a locally $\delta$-compact $T_{3}$ space $(X,\tau)$.\\
let $y\in Y$. Then $y\in X$.\\
Since $X$ is locally $\delta$- compact, then $\exists$ a neighbourhood $V$ of $y$ in $(X,\tau)$ which is $\delta$-compact in $X$.\\
Since $Y$ is closed of a regular space $X$, $Y$ is a $\delta$-closed set by the Theorem \ref{th8}. Again since $V$ is $\delta$-compact subset of a $T_{3}$ space $X$, $V$ is $\delta$-closed by the Theorem \ref{th15}. So, $V\cap Y$is $\delta$-closed subset of a $\delta$-compact $V$. Thus by the Theorem \ref{th7}, $V\cap Y$ is $\delta$-compact. So, $V\cap Y$ is a $\delta$-compact neighbour of  $y$ in $Y$. 
Hence, by arbitrariness of $y$, it proves that $(Y,\tau_{Y})$ is locally $\delta$-compact.
\end{proof}

\section{Conclusion}
In this paper, at first we establish some basic results on $\delta$-Compact sets with the help of regular open sets, regular closed sets, $\delta$-open sets, $\delta$-closed sets. Then we have examined the behavior of $\delta$-Compact sets under a continuous mapping. Lastly, we define locally $\delta$-compact space and by these, we prove that every closed subspace of a locally $\delta$-compact $T_{3}$ space is locally $\delta$-compact.\\
So, one can try to establish the Alexender Sub-Base Theorem with respect to $\delta$-Compact space.

\end{document}